\documentclass[12pt,a4paper,dvips]{article}
\usepackage{inputenc, amssymb, amsmath, amsthm, srcltx}
\usepackage[russian,english]{babel}
\usepackage{color}
\usepackage{graphicx}
\makeatletter
\thm@style{plain}

\def\th@plain{%
  \thm@headfont{\bfseries}%
  \itshape 
  \thm@notefont{\rm}%
}
\def\thm@indent{\hspace*{\parindent}}
\makeatother

\def\({\left(}
\def\){\right)}

\newcommand{\be}{\begin{equation}}
\newcommand{\ee}{\end{equation}}

\let\epsilon\varepsilon
\let\phi\varphi
\let\le\leqslant
\let\ge\geqslant

\def\arg{\mathop{\ensuremath{\text{\textup{arg}}}}}

\def\Re{\mathop{\ensuremath{\text{\textup{Re}}}}}
\def\Im{\mathop{\ensuremath{\text{\textup{Im}}}}}

\newtheorem{theorem}{Theorem}
\newtheorem{lemma}{Lemma}

\newcommand{\beq}{\begin{equation}}
\newcommand{\eeq}{\end{equation}}

{Definition}
{Algorithm}

\begin{document}

\centerline{\bf\uppercase{A small improvement in the small gaps}}
\centerline{\bf\uppercase{between consecutive zeros of the Riemann zeta-function}\footnote[1]{%
}}

\bigskip

\medskip

\centerline{\sc Sergei~Preobrazhenski\u i}

\bigskip

\bigskip

\hbox to \textwidth{\hfil\parbox{0.9\textwidth}{%
\small {\sc Abstract.}
Feng and Wu introduced a new general coefficient sequence into
Montgomery and Odlyzko's method for exhibiting irregularity
in the gaps between consecutive zeros of $\zeta(s)$
assuming the Riemann Hypothesis.
They used a special case of their sequence to improve upon earlier results on the gaps.
In this paper we consider an equivalent form of the general sequence of Feng and Wu,
and introduce a somewhat less general sequence $\{a_n\}$
for which we write the Montgomery--Odlyzko expressions explicitly.
As an application, we give the following slight improvement of Feng and Wu's result:
infinitely often consecutive non-trivial
zeros of the Riemann zeta-function differ by at most $0{.}515396$
times the average spacing.}\hfil}

\bigskip

\bigskip
\textbf{Contents}\hfil

\bigskip

1. Introduction

\smallskip

2. Lemmas

\smallskip

3. Proof of Theorem 1

\smallskip

References

\bigskip

\bigskip

\bigskip


{\bf 1. Introduction.} It is well known that the Riemann zeta-function $\zeta(s)$
has infinitely many nontrivial zeros $s=\rho=\beta+i\gamma$,
and all of them are in the critical strip
$0<\Re s=\sigma<1$, $-\infty<\Im s=t<\infty$.

If $N(T)$ denotes the number of zeros $\rho=\beta+i\gamma$
($\beta$ and $\gamma$ real), for which $0<\gamma\le T$,
then
\[
N(T)=\frac T{2\pi}\log\left(\frac T{2\pi}\right)
-\frac T{2\pi}+\frac78+S(T)+O\left(\frac1T\right),
\]
with
\[
S(T)=\frac1{\pi}\arg\zeta\left(\frac12+iT\right)
\]
and
\[
S(T)=O(\log T).
\]
This is the Riemann--von~Mangoldt formula for $N(T)$.
Hence, if we let $0<\gamma\le\gamma'$ denote consecutive ordinates
of non-trivial zeros of $\zeta(s)$, the average size of $\gamma'-\gamma$
is $\gamma/N(\gamma)\sim2\pi/\log\gamma$. Let
\[
\lambda=\limsup\limits_{\gamma>0}(\gamma'-\gamma)\frac{\log\gamma}{2\pi}
\]
and
\[
\mu=\liminf\limits_{\gamma>0}(\gamma'-\gamma)\frac{\log\gamma}{2\pi}.
\]
We note that $\mu\le1\le\lambda$ and it is expected that $\mu=0$ and $\lambda=+\infty$.

Let $N_0(T)$ be the number of zeros of $\zeta\left(\frac12+it\right)$ when $0<t\le T$,
each zero counted with multiplicity.
The Riemann hypothesis is the conjecture that $N_0(T)=N(T)$.

In this note we prove the following theorem.
\begin{theorem}\label{PreobSmallGaps}
Assume the Riemann Hypothesis. Then we have
\[
\mu<0{.}515396.
\]
\end{theorem}

We briefly describe the history of the problem, focusing mainly on $\mu$.
\begin{itemize}
\item \cite{PreobSel46}: in 1946 Selberg remarked that $\mu<1<\lambda$ unconditionally.
\end{itemize}

Now suppose that $T$ is a large real number and $K=T(\log T)^{-2}$. Let
\be\label{Preobhc}
h(c)=c-\frac{\textrm{Re}\left(\sum_{nk\le K}a_k\overline{a_{nk}}g_c(n)\Lambda(n)n^{-1/2}\right)}
{\sum_{k\le K}\left\vert a_k\right\vert^2}
\ee
where
\[
g_c(n)=\frac{2\sin\left(\pi c\frac{\log n}{\log T}\right)}{\pi\log n}
\]
and $\Lambda$ is the von Mangoldt's function.

In the following results, the truth of the Riemann Hypothesis is assumed.
\begin{itemize}
\item \cite{PreobMO81}: in 1981 by an argument using the Guinand--Weil explicit formula, Montgomery and Odlyzko
showed that if $h(c)<1$ for some choice of $c$ and $\{a_n\}$, then $\lambda\ge c$, and if $h(c)>1$ for some choice of $c$
and $\{a_n\}$, then $\mu\le c$. They used the coefficients
\[
a_k=\frac1{k^{\frac12}}f\left(\frac{\log k}{\log K}\right)\quad\textrm{ and }\quad
a_k=\frac{\lambda(k)}{k^{\frac12}}f\left(\frac{\log k}{\log K}\right)
\]
where $f$ is a continuous function of bounded variation, and $\lambda(k)$ is the Liouville function.
With this choice of the coefficients they obtained
$\lambda>1{.}9799$ and $\mu<0{.}5179$ by optimizing over the functions $f$.

\item \cite{PreobCGG84}: in 1984 Conrey, Ghosh \& Gonek chose the coefficients
\[
a_k=\frac{d_r(k)}{\sqrt{k}}\quad\text{ and }\quad a_k=\frac{\lambda(k) d_r(k)}{\sqrt{k}}
\]
where $d_r(k)$ is a multiplicative function defined on integral powers of a prime $p$ by
\[
d_r(p^k)=\frac{\Gamma(k+r)}{\Gamma(r) k!}.
\]
The choice $r=1{.}1$ with the latter $a_k$ yields $\mu<0{.}5172$ and the choice $r=2{.}2$
with the former $a_k$ yields $\lambda>2{.}337$.

\item \cite{PreobHall05}: in 2005, by making use of the Wirtinger's inequality
and the asymptotic formulae for the fourth mixed moments of the zeta-function and its derivative,
R.~R.~Hall proved that $\lambda>2{.}6306$.

\item \cite{PreobBMN10}: in 2010, Bui, Milinovich \& Ng
considered the coefficients of the form
\[
a_k=\frac{d_r(k)}{\sqrt{k}} f\big(\tfrac{\log K/k}{\log K}\big)\quad\text{ and }\quad
a_k=\frac{\lambda(k) d_r(k)}{\sqrt{k}}f\big(\tfrac{\log K/k}{\log K}\big)
\]
for a polynomial $f$ and obtained $\lambda>2{.}69$ and $\mu<0{.}5155$.

\item \cite{PreobFW12}: in 2012, Feng \& Wu introduced the coefficient
\[
\begin{split}
a_k=&\frac{d_r(k)}{k^{\frac12}}\left(f_1\left(\frac{\log K/k}{\log
K}\right)+f_2\left(\frac{\log K/k}{\log K}\right)\sum_{p_1p_2\mid k}\frac{\log
p_1\log p_2}{\log^2K}\right.\\
&+f_3\left(\frac{\log K/k}{\log K}\right)\sum_{p_1p_2p_3\mid k}\frac{\log
p_1\log p_2\log p_3}{\log^3K}+\cdots\\
&\left.+f_I\left(\frac{\log K/k}{\log K}\right)\sum_{p_1p_2\cdots p_I\mid k}\frac{\log
p_1\log p_2\cdots\log p_I}{\log^IK}\right),
\end{split}
\]
for any integer $I\ge2$. Using $I=2$ they obtained $\lambda>2{.}7327$ and $\mu<0{.}5154$,
or, to higher precision, $\mu<0{.}515398$.

\end{itemize}

We remark that the coefficient of Feng \& Wu is equivalent to the coefficient
\[
\begin{split}
a_k=&\frac{d_r(k)}{k^{\frac12}}\left(f_1\left(\frac{\log K/k}{\log
K}\right)+f_2\left(\frac{\log K/k}{\log K}\right)\sum_{p_1\mid k}\frac{\log^2
p_1}{\log^2K}\right.\\
&+f_3\left(\frac{\log K/k}{\log K}\right)\sum_{p_1\mid k}\frac{\log^3
p_1}{\log^3K}+\cdots\\
&\left.+f_I\left(\frac{\log K/k}{\log K}\right)\sum_{p_1\mid k}\frac{\log^I
p_1}{\log^IK}\right),
\end{split}
\]
for which the calculations are simpler.

To prove Theorem~\ref{PreobSmallGaps}, we choose the coefficients
\[
a_k=\frac{\lambda(k)d_r(k)}{k^{\frac12}}f_1\left(\frac{\log K/k}{\log K}\right)
+\frac{\lambda(k)d_r(k)}{k^{\frac12}}\sum_{p\mid k}P\left(\frac{\log p}{\log K}\right)
\tilde{f}_1\left(\frac{\log K/k}{\log K}\right),
\]
where $f_1$, $\tilde{f}_1$, $P$ are some polynomials to be chosen later.
The $a_k$ given above are less general than the coefficients of Feng and Wu,
but they are simpler, so we are able to write the Montgomery--Odlyzko expressions
for our sequence explicitly.

\textbf{2. Lemmas.}

\begin{lemma}[Mertens Theorem]
\label{PreobMertens}
\[
\sum_{p\le y}\frac{\log p}{p}=\log y+O(1).
\]
\end{lemma}

\begin{lemma}[See Levinson~\cite{PreobLev74}]
\label{PreobLev}
\[
\sum_{p\mid j}\frac{\log p}{p}=O(\log\log j).
\]
\end{lemma}

\begin{lemma}
\label{Preobdr2}
For fixed $r\ge1$,
\[
\sum_{k\le x}\frac{d_r(k)^2}{k}=A_r(\log x)^{r^2}+O\left((\log T)^{r^2-1}\right)
\]
uniformly for $x\le T$.
\end{lemma}

\begin{lemma}
\label{Preobdimreduct}
Let $a_i$ be integer for $1\le i\le m$, $D>1$ and $f$
is a continuous function. Then
\[
\begin{split}
&\int_{1}^{D}\frac{\log ^{a_1-1}x_1}{x_1}dx_1\int_{1}^{\frac
{D}{x_1}}\frac{\log ^{a_2-1}x_2}{x_2}dx_2\cdots\int_1^{\frac D{x_1x_2\cdots x_m}}\frac{f(x_1x_2\cdots x_mx)}xdx\\
=&\frac{\prod_{i=1}^m(a_i-1)!}{(\sum_{i=1}^ma_i)!}\int_1^D\frac{f(x)\log
^{\sum_{i=1}^ma_i}x}xdx.
\end{split}
\]
\end{lemma}

\begin{lemma}
\label{PreobFengWu}
Let $a_i$ be integer for $1\le i\le m$, and $g$
is a polynomial. Then we have
\[
\begin{split}
&\sum_{k\le K}\frac{d_r(k)^2}kg\left(\frac{\log K/k}{\log
K}\right)\\
=&A_rr^2\int_1^Kg\left(\frac{\log K/x}{\log K}\right)(\log x)^{r^2-1}\frac{dx}{x}+O\left((\log K)^{r^2-1}\right)
\end{split}
\]
and
\[
\begin{split}
&\sum_{p_1p_2\cdots p_m\le K}\prod_{i=1}^{m}\frac{\log^{a_i}p_i}{p_i}\mu^2(p_1p_2\cdots p_m)
\sum_{k_0\le K/(p_1p_2\dots p_m)}\frac{d_r(k_0)^2}{k_0}g\left(\frac{\log K/(p_1p_2\dots p_mk_0)}{\log
K}\right)\\
=&\left(1+O\left(\log^{-1}K\right)\right)\\
&\times A_rr^2
\int_1^K\log^{a_1-1}x_1\frac{dx_1}{x_1}\int_1^{\frac{K}{x_1}}\log^{a_2-1}x_2\frac{dx_2}{x_2}
\cdots\int_1^{\frac{K}{x_1x_2\cdots x_{m-1}}}\log^{a_m-1}x_m\frac{dx_m}{x_m}\\
&\times\int_1^{\frac{K}{x_1x_2\cdots x_m}}g\left(\frac{\log K/(x_1x_2\cdots x_mx)}{\log K}\right)(\log x)^{r^2-1}\frac{dx}{x}.
\end{split}
\]
\end{lemma}

For the proof of Lemma~\ref{PreobFengWu} using Lemmas~\ref{PreobMertens}--\ref{Preobdr2}, see~\cite{PreobFW12}.

\textbf{3. Proof of Theorem~\ref{PreobSmallGaps}.}
To give an upper bound for $\mu$, we evaluate $h(c)$ in~\eqref{Preobhc} with the coefficients
\[
a_k=\frac{\lambda(k)d_r(k)}{k^{\frac12}}f_1\left(\frac{\log K/k}{\log K}\right)
+\frac{\lambda(k)d_r(k)}{k^{\frac12}}\sum_{p\mid k}P\left(\frac{\log p}{\log K}\right)
\tilde{f}_1\left(\frac{\log K/k}{\log K}\right),
\]
where $r\ge1$ and $f_1$, $\tilde{f}_1$, $P$ are polynomials.

First, we evaluate the denominator in the ratio in the definition of $h(c)$.
\[
\begin{split}
\sum_{k\le K}\left\vert a_k\right\vert^2=&\sum_{k\le K}\frac{d_r(k)^2}kf_1\left(\frac{\log K/k}{\log
K}\right)^2\\
&+2\sum_{k\le K}\frac{d_r(k)^2}kf_1\left(\frac{\log K/k}{\log K}\right)\tilde{f}_1\left(\frac{\log K/k}{\log K}\right)
\sum_{p\mid k}P\left(\frac{\log p}{\log K}\right)\\
&+\sum_{k\le K}\frac{d_r(k)^2}k\tilde{f}_1\left(\frac{\log K/k}{\log
K}\right)^2\sum_{p\mid k}P\left(\frac{\log p}{\log K}\right)
\sum_{q\mid k}P\left(\frac{\log q}{\log K}\right)\\
=&\tilde{D}_1+\tilde{D}_2+\tilde{D}_3.
\end{split}
\]
Using Lemma~\ref{PreobFengWu} and recalling that $K=T(\log T)^{-2}$, we have
\be
\label{PreobD1}
\begin{split}
\tilde{D}_1=&A_rr^2\int_1^Kf_1\left(\frac{\log K/x}{\log K}\right)^2(\log x)^{r^2-1}\frac{dx}{x}+O\left((\log T)^{r^2-1}\right)\\
=&A_rr^2(\log K)^{r^2}\int_0^1(1-u)^{r^2-1}f_1(u)^2du+O\left((\log T)^{r^2-1}\right)\\
=&A_rr^2(\log T)^{r^2}\int_0^1(1-u)^{r^2-1}f_1(u)^2du+O\left((\log T)^{r^2-1+\epsilon}\right),
\end{split}
\ee
where $\epsilon>0$ is arbitrarily small and the constant in the $O$-term depends on $r$, $\epsilon$ and $f_1$.
By Lemma~\ref{PreobFengWu} we obtain that
\[
\begin{split}
\tilde{D}_2=&\frac{2A_rr^4}{\log K}\int_1^K\frac{P_1(\log x_1)}{x_1}\int_1^{\frac{K}{x_1}}f_1\left(\frac{\log K/x_1x}{\log K}\right)\\
&\times\tilde{f}_1\left(\frac{\log K/x_1x}{\log K}\right)(\log x)^{r^2-1}\frac{dx}{x}dx_1+O\left((\log T)^{r^2-1+\epsilon}\right),
\end{split}
\]
where $P_1(y)=\frac{P(y)}y$.
By the variable changes $u=1-\frac{\log x_1}{\log K}$, $v=1-\frac{\log x_1x}{\log K}$, we have
\be
\label{PreobD2}
\begin{split}
\tilde{D}_2=&2A_rr^4(\log K)^{r^2}\int_0^1P_1(1-u)\int_0^u(u-v)^{r^2-1}f_1(v)\tilde{f}_1(v)dvdu+O\left((\log T)^{r^2-1+\epsilon}\right)\\
=&2A_rr^4(\log T)^{r^2}\int_0^1P_1(1-u)\int_0^u(u-v)^{r^2-1}f_1(v)\tilde{f}_1(v)dvdu+O\left((\log T)^{r^2-1+\epsilon}\right),
\end{split}
\ee
where the constant in the $O$-term depends on $r$, $\epsilon$ and $f_1$, $\tilde{f}_1$.

We have
\[
\begin{split}
\tilde{D}_3=&\sum_{k\le K}\frac{d_r(k)^2}k\tilde{f}_1\left(\frac{\log K/k}{\log K}\right)^2
\sum_{p_1p_2\mid k}\mu^2(p_1p_2)P\left(\frac{\log p_1}{\log K}\right)P\left(\frac{\log p_2}{\log K}\right)\\
&+\sum_{k\le K}\frac{d_r(k)^2}k\tilde{f}_1\left(\frac{\log K/k}{\log
K}\right)^2\sum_{p\mid k}P^2\left(\frac{\log p}{\log K}\right)\\
=&\tilde{D}_{31}+\tilde{D}_{32}.
\end{split}
\]
Again by Lemma~\ref{PreobFengWu} we obtain that
\[
\begin{split}
\tilde{D}_{31}=&\frac{A_rr^6}{\log^2K}\int_1^K\frac{P_1(\log x_1)}{x_1}\int_1^{\frac{K}{x_1}}\frac{P_1(\log x_2)}{x_2}
\int_1^{\frac{K}{x_1x_2}}\tilde{f}_1\left(\frac{\log K/x_1x_2x}{\log K}\right)^2\\
&\times(\log x)^{r^2-1}\frac{dx}{x}dx_2dx_1+O\left((\log T)^{r^2-1+\epsilon}\right),
\end{split}
\]
where $P_1(y)=\frac{P(y)}y$.
We remark that by Lemma~\ref{Preobdimreduct} we can reduce the number of the repeated integrations in the above expression.
By the change of variables $u=1-\frac{\log x_1}{\log K}$, $v=1-\frac{\log x_2}{\log K}$, $w=1-\frac{\log x_1x_2x}{\log K}$,
\be
\label{PreobD31}
\begin{split}
\tilde{D}_{31}=&A_rr^6(\log K)^{r^2}\int_0^1P_1(1-u)\int_{1-u}^1P_1(1-v)\int_0^{u+v-1}(u+v-w-1)^{r^2-1}\tilde{f}_1(w)^2dw\,dv\,du\\
&+O\left((\log T)^{r^2-1+\epsilon}\right)\\
=&A_rr^6(\log T)^{r^2}\int_0^1P_1(1-u)\int_{1-u}^1P_1(1-v)\int_0^{u+v-1}(u+v-w-1)^{r^2-1}\tilde{f}_1(w)^2dw\,dv\,du\\
&+O\left((\log T)^{r^2-1+\epsilon}\right),
\end{split}
\ee
where the constant in the $O$-term depends on $r$, $\epsilon$ and $\tilde{f}_1$.
Similarly,
\be
\label{PreobD32}
\begin{split}
\tilde{D}_{32}=&A_rr^4(\log K)^{r^2}\int_0^1P_2(1-u)\int_0^u(u-v)^{r^2-1}\tilde{f}_1(v)^2dv\,du+O\left((\log T)^{r^2-1+\epsilon}\right)\\
=&A_rr^4(\log T)^{r^2}\int_0^1P_2(1-u)\int_0^u(u-v)^{r^2-1}\tilde{f}_1(v)^2dv\,du+O\left((\log T)^{r^2-1+\epsilon}\right),
\end{split}
\ee
where $P_2(y)=\frac{P(y)^2}y$ and the constant in the $O$-term depends on $r$, $\epsilon$ and $\tilde{f}_1$.

We now proceed to evaluation of the numerator in the ratio in~\eqref{Preobhc}.
If we let
\[
N(c)=\sum_{nk\le K}a_ka_{nk}g_c(n)\Lambda(n)n^{-1/2},
\]
then
\[
\begin{split}
N(c)=&\frac2{\pi}\sum_{nk\le
K}\frac{\lambda(k)d_r(k)\lambda(nk)d_r(nk)\Lambda(n)}{kn\log n}\sin\left(\pi c\frac{\log
n}{\log T}\right)\times\left(f_1\left(\frac{\log K/k}{\log K}\right)f_1\left(\frac{\log K/nk}{\log
K}\right)\right.\\
&+f_1\left(\frac{\log K/nk}{\log K}\right)\tilde{f}_1\left(\frac{\log K/k}{\log
K}\right)\sum_{p_1\mid k}P\left(\frac{\log p_1}{\log K}\right)\\
&+f_1\left(\frac{\log K/k}{\log K}\right)\tilde{f}_1\left(\frac{\log K/nk}{\log
K}\right)\sum_{p_1\mid nk}P\left(\frac{\log p_1}{\log K}\right)\\
&\left.{}+\tilde{f}_1\left(\frac{\log K/k}{\log
K}\right)\tilde{f}_1\left(\frac{\log K/nk}{\log K}\right)\sum_{p_1\mid k}P\left(\frac{\log p_1}{\log K}\right)
\sum_{q_1\mid nk}P\left(\frac{\log q_1}{\log K}\right)\right),
\end{split}
\]
so we can write
\[
N(c)=N_1+N_2+N_3+N_4.
\]
Using the distribution of $\Lambda(n)$, we obtain
\[
\begin{split}
N_1=&-\frac2{\pi}\sum_{pk\le K}\frac{d_r(k)d_r(pk)}{kp}\sin\left(\pi c\frac{\log
p}{\log T}\right)f_1\left(\frac{\log K/k}{\log K}\right)f_1\left(\frac{\log K/pk}{\log
K}\right)+O\left((\log T)^{r^2-1}\right)\\
=&-\frac{2r}{\pi}\sum_{p\le K}\frac{\sin\left(\pi c\frac{\log
p}{\log T}\right)}{p}\sum_{k\le K/p}\frac{d_r(k)^2}{k}f_1\left(\frac{\log K/k}{\log K}\right)f_1\left(\frac{\log K/pk}{\log
K}\right)+O\left((\log T)^{r^2-1}\right).
\end{split}
\]
By Lemma~\ref{PreobFengWu} we have
\[
\begin{split}
N_1=&-\frac{2A_rr^3}{\pi}\sum_{p\le K}\frac{\sin\left(\pi c\frac{\log p}{\log T}\right)}{p}
\int_1^{\frac K{p}}f_1\left(\frac{\log K/x}{\log K}\right)f_1\left(\frac{\log K/px}{\log K}\right)(\log x)^{r^2-1}\frac{dx}{x}\\
&+O\left((\log T)^{r^2-1}\right).
\end{split}
\]
From Lemma~\ref{PreobMertens} and Abel's summation,
\[
\begin{split}
N_1=&-\frac{2A_rr^3}{\pi}\int_1^K\frac{\sin\left(\pi c\frac{\log x_1}{\log T}\right)}{x_1\log x_1}\int_1^{\frac K{x_1}}
f_1\left(\frac{\log K/x}{\log K}\right)f_1\left(\frac{\log K/xx_1}{\log K}\right)(\log x)^{r^2-1}\frac{dx}{x}dx_1\\
&+O\left((\log T)^{r^2-1}\right).
\end{split}
\]
Interchanging the order of integration and the names of the variables $x$ and $x_1$, we find
\[
\begin{split}
N_1=&-\frac{2A_rr^3}{\pi}\int_1^Kf_1\left(\frac{\log K/x_1}{\log K}\right)\frac{(\log x_1)^{r^2-1}}{x_1}\int_1^{\frac K{x_1}}
\frac{\sin\left(\pi c\frac{\log x}{\log T}\right)}{\log x}f_1\left(\frac{\log K/xx_1}{\log K}\right)\frac{dx}{x}dx_1\\
&+O\left((\log T)^{r^2-1}\right).
\end{split}
\]
Let $u=1-\frac{\log x_1}{\log K}$, $v=\frac{\log x}{\log K}$. Then
\be
\label{PreobN1}
\begin{split}
N_1=&-\frac{2A_rr^3}{\pi}(\log K)^{r^2}\int_0^1(1-u)^{r^2-1}f_1(u)\int_0^u\frac{\sin\left(\pi
cv\frac{\log K}{\log T}\right)}{v}f_1(u-v)\,dv\,du\\
&+O\left((\log T)^{r^2-1}\right)\\
=&-\frac{2A_rr^3}{\pi}(\log T)^{r^2}\int_0^1(1-u)^{r^2-1}f_1(u)\int_0^u\frac{\sin(\pi
cv)}{v}f_1(u-v)\,dv\,du\\
&+O\left((\log T)^{r^2-1+\epsilon}\right),
\end{split}
\ee
where the constant in the $O$-term depends on $r$, $\epsilon$ and $f_1$.

In $N_2$ we can replace the product of the summation variables $nk$ by $pp_1k_0$ to get
\[
\begin{split}
N_2=&-\frac{2r^3}{\pi}\sum_{p_1\le K}\frac{P\left(\frac{\log p_1}{\log K}\right)}{p_1}\sum_{pk_0\le K/p_1}\frac{\sin\left(\pi c\frac{\log
p}{\log T}\right)d_r(k_0)^2}{pk_0}\\
&\times f_1\left(\frac{\log K/(pp_1k_0)}{\log K}\right)\tilde{f}_1\left(\frac{\log K/(p_1k_0)}{\log
K}\right)+O\left((\log T)^{r^2-1}\right).
\end{split}
\]
The inner sum $\sum_{pk_0\le K/p_1}$ in the expression above
is the sum $\sum_{k_0\le K/p_1}\sum_{p\le K/(p_1k_0)}$.
As in the calculation of $N_1$, we can show that this double sum is
\[
\begin{split}
&A_rr^2\int_1^{\frac K{p_1}}\tilde{f}_1\left(\frac{\log K/(p_1x_2)}{\log K}\right)\frac{(\log x_2)^{r^2-1}}{x_2}\\
&\times\int_1^{\frac K{p_1x_2}}
\frac{\sin\left(\pi c\frac{\log x}{\log T}\right)}{\log x}f_1\left(\frac{\log K/(p_1xx_2)}{\log K}\right)\frac{dx}{x}dx_2
+O\left((\log T)^{r^2-1}\right).
\end{split}
\]
By Lemma~\ref{PreobMertens} we obtain
\[
\begin{split}
N_2=&-\frac{2A_rr^5}{\pi(\log K)}\int_1^K\frac{P_1\left(\frac{\log x_1}{\log K}\right)}{x_1}
\int_1^{\frac K{x_1}}\tilde{f}_1\left(\frac{\log K/(x_1x_2)}{\log K}\right)\frac{(\log x_2)^{r^2-1}}{x_2}\\
&\times\int_1^{\frac K{x_1x_2}}\frac{\sin\left(\pi c\frac{\log x}{\log T}\right)}{\log x}f_1\left(\frac{\log K/(xx_1x_2)}{\log K}\right)
\frac{dx}{x}dx_2\,dx_1+O\left((\log T)^{r^2-1}\right).
\end{split}
\]
Making the variable changes $u=1-\frac{\log x_1}{\log K}$, $v=1-\frac{\log x_1x_2}{\log K}$, $w=\frac{\log x}{\log K}$, we get
\be
\label{PreobN2}
\begin{split}
N_2=&-\frac{2A_rr^5}{\pi}(\log T)^{r^2}\int_0^1P_1(1-u)\int_0^u(u-v)^{r^2-1}\tilde{f}_1(v)\int_0^v\frac{\sin(\pi cw)}{w}f_1(v-w)\,dw\,dv\,du\\
&+O\left((\log T)^{r^2-1+\epsilon}\right),
\end{split}
\ee
where the constant in the $O$-term depends on $r$, $\epsilon$ and $f_1$, $\tilde{f}_1$.

As in $N_1$ and $N_2$, the terms with $n=p$ for the primes $p$ give the main contribution to $N_3$:
\[
\begin{split}
N_3=&-\frac2{\pi}\sum_{pk\le K}\sin\left(\pi c\frac{\log
p}{\log T}\right)\frac{d_r(k)d_r(kp)}{kp}f_1\left(\frac{\log K/k}{\log K}\right)\tilde{f}_1\left(\frac{\log
K/(pk)}{\log K}\right)\\
&\times\sum_{p_1\mid pk}P\left(\frac{\log p_1}{\log K}\right)+O\left((\log T)^{r^2-1}\right).
\end{split}
\]
For $(p,k)=1$ it follows that
\be
\label{PreobN3decom}
\sum_{p_1\mid pk}P\left(\frac{\log p_1}{\log K}\right)=\sum_{p_1\mid k}P\left(\frac{\log p_1}{\log K}\right)+P\left(\frac{\log p}{\log K}\right).
\ee
Since the contribution of the terms with $(p,k)\neq1$ in $N_3$ is $O\left((\log T)^{r^2-1}\right)$,
then, according to decomposition~\eqref{PreobN3decom}, we can write
\[
N_3=N_{31}+N_{32}+O\left((\log T)^{r^2-1}\right),
\]
where
\[
\begin{split}
N_{31}=&-\frac{2r^3}{\pi}\sum_{p_1\le K}\frac{P\left(\frac{\log p_1}{\log K}\right)}{p_1}
\sum_{pk_0\le K/p_1}\frac{\sin\left(\pi c\frac{\log
p}{\log T}\right)d_r(k_0)^2}{pk_0}\\
&\times\tilde{f}_1\left(\frac{\log K/(pp_1k_0)}{\log K}\right)f_1\left(\frac{\log K/(p_1k_0)}{\log
K}\right)+O\left((\log T)^{r^2-1}\right)
\end{split}
\]
and
\[
\begin{split}
N_{32}=&-\frac{2r}{\pi}\sum_{pk\le K}\frac{\sin\left(\pi c\frac{\log
p}{\log T}\right)d_r(k)^2P\left(\frac{\log p}{\log K}\right)}{pk}\tilde{f}_1\left(\frac{\log K/(pk)}{\log K}\right)\\
&\times f_1\left(\frac{\log K/k}{\log
K}\right)+O\left((\log T)^{r^2-1}\right).
\end{split}
\]
As in the calculation of $N_2$ we get
\[
\begin{split}
N_{31}=&-\frac{2A_rr^5}{\pi}(\log T)^{r^2}\int_0^1P_1(1-u)\int_0^u(u-v)^{r^2-1}f_1(v)\\
&\times\int_0^v\frac{\sin(\pi cw)}{w}\tilde{f}_1(v-w)\,dw\,dv\,du+O\left((\log T)^{r^2-1+\epsilon}\right),
\end{split}
\]
and as in the calculation of $N_1$,
\[
\begin{split}
N_{32}=&-\frac{2A_rr^3}{\pi}(\log T)^{r^2}\int_0^1(1-u)^{r^2-1}f_1(u)\int_0^u\sin(\pi cv)P_1(v)\tilde{f}_1(u-v)\,dv\,du\\
&+O\left((\log T)^{r^2-1+\epsilon}\right),
\end{split}
\]
where $P_1(y)=\frac{P(y)}y$.
Thus,
\be
\label{PreobN3}
\begin{split}
N_3=&-\frac{2A_rr^5}{\pi}(\log T)^{r^2}\int_0^1P_1(1-u)\int_0^u(u-v)^{r^2-1}f_1(v)\\
&\times\int_0^v\frac{\sin(\pi cw)}{w}\tilde{f}_1(v-w)\,dw\,dv\,du\\
&-\frac{2A_rr^3}{\pi}(\log T)^{r^2}\int_0^1(1-u)^{r^2-1}f_1(u)\int_0^u\sin(\pi cv)P_1(v)\tilde{f}_1(u-v)\,dv\,du\\
&+O\left((\log T)^{r^2-1+\epsilon}\right),
\end{split}
\ee
where the constant in the $O$-term depends on $r$, $\epsilon$ and $f_1$, $\tilde{f}_1$, $P$.

Again, in the sum defining $N_4$ we can replace the integers $n\ge2$ with the primes $p$:
\[
\begin{split}
N_4=&-\frac2{\pi}\sum_{pk\le K}\sin\left(\pi c\frac{\log
p}{\log T}\right)\frac{d_r(k)d_r(kp)}{kp}\tilde{f}_1\left(\frac{\log K/k}{\log K}\right)\tilde{f}_1\left(\frac{\log K/(pk)}{\log K}\right)\\
&\times\sum_{p_1\mid k}P\left(\frac{\log p_1}{\log K}\right)\sum_{q_1\mid pk}P\left(\frac{\log q_1}{\log K}\right)
+O\left((\log T)^{r^2-1}\right).
\end{split}
\]
For the two innermost sums, if $(k,p)=1$, we have
\[
\begin{split}
&\sum_{p_1\mid k}P\left(\frac{\log p_1}{\log K}\right)\sum_{q_1\mid pk}P\left(\frac{\log q_1}{\log K}\right)\\
=&\sum_{p_1q_1\mid k}\mu^2(p_1q_1)P\left(\frac{\log p_1}{\log K}\right)P\left(\frac{\log q_1}{\log K}\right)\\
&+\sum_{p_1\mid k}P^2\left(\frac{\log p_1}{\log K}\right)+P\left(\frac{\log p}{\log K}\right)\sum_{p_1\mid k}P\left(\frac{\log p_1}{\log K}\right).
\end{split}
\]
According to this decomposition, we write
\[
N_4=N_{41}+N_{42}+N_{43}.
\]
As before, by Lemma~\ref{PreobFengWu} we find
\[
\begin{split}
N_{41}=&-\frac{2A_rr^7}{\pi(\log K)^2}\int_1^K\frac{P_1\left(\frac{\log x_1}{\log K}\right)}{x_1}
\int_1^{\frac K{x_1}}\frac{P_1\left(\frac{\log x_2}{\log K}\right)}{x_2}
\int_1^{\frac K{x_1x_2}}\tilde{f}_1\left(\frac{\log K/(x_1x_2x_3)}{\log K}\right)\frac{(\log x_3)^{r^2-1}}{x_3}\\
&\times\int_1^{\frac K{x_1x_2x_3}}\frac{\sin\left(\pi c\frac{\log x}{\log T}\right)}{\log x}\tilde{f}_1\left(\frac{\log K/(xx_1x_2x_3)}{\log K}\right)
\frac{dx}{x}dx_3\,dx_2\,dx_1+O\left((\log T)^{r^2-1}\right).
\end{split}
\]
Making the variable changes $u=1-\frac{\log x_1}{\log K}$, $v=1-\frac{\log x_2}{\log K}$, $w=1-\frac{\log x_1x_2x_3}{\log K}$, $z=\frac{\log x}{\log K}$, we get
\be
\label{PreobN41}
\begin{split}
N_{41}=&-\frac{2A_rr^7}{\pi}(\log T)^{r^2}\int_0^1P_1(1-u)\int_{1-u}^1P_1(1-v)\\
&\times\int_0^{u+v-1}(u+v-w-1)^{r^2-1}\tilde{f}_1(w)\int_0^w\frac{\sin(\pi cz)}{z}\tilde{f}_1(w-z)\,dz\,dw\,dv\,du\\
&+O\left((\log T)^{r^2-1+\epsilon}\right),
\end{split}
\ee
where the constant in the $O$-term depends on $r$, $\epsilon$ and $\tilde{f}_1$, $P$.

Next,
\[
\begin{split}
N_{42}=&-\frac{2A_rr^5}{\pi(\log K)}\int_1^K\frac{P_2\left(\frac{\log x_1}{\log K}\right)}{x_1}
\int_1^{\frac K{x_1}}\tilde{f}_1\left(\frac{\log K/(x_1x_2)}{\log K}\right)\frac{(\log x_2)^{r^2-1}}{x_2}\\
&\times\int_1^{\frac K{x_1x_2}}\frac{\sin\left(\pi c\frac{\log x}{\log T}\right)}{\log x}\tilde{f}_1\left(\frac{\log K/(xx_1x_2)}{\log K}\right)
\frac{dx}{x}dx_2\,dx_1+O\left((\log T)^{r^2-1}\right).
\end{split}
\]
By the variable changes $u=1-\frac{\log x_1}{\log K}$, $v=1-\frac{\log x_1x_2}{\log K}$, $w=\frac{\log x}{\log K}$, we get
\be
\label{PreobN42}
\begin{split}
N_{42}=&-\frac{2A_rr^5}{\pi}(\log T)^{r^2}\int_0^1P_2(1-u)\int_0^u(u-v)^{r^2-1}\tilde{f}_1(v)\int_0^v\frac{\sin(\pi cw)}{w}\tilde{f}_1(v-w)\,dw\,dv\,du\\
&+O\left((\log T)^{r^2-1+\epsilon}\right),
\end{split}
\ee
where the constant in the $O$-term depends on $r$, $\epsilon$ and $\tilde{f}_1$, $P$.

Finally,
\[
\begin{split}
N_{43}=&-\frac{2A_rr^5}{\pi(\log K)^2}\int_1^K\frac{P_1\left(\frac{\log x_1}{\log K}\right)}{x_1}
\int_1^{\frac K{x_1}}\tilde{f}_1\left(\frac{\log K/(x_1x_2)}{\log K}\right)\frac{(\log x_2)^{r^2-1}}{x_2}\\
&\times\int_1^{\frac K{x_1x_2}}\sin\left(\pi c\frac{\log x}{\log T}\right)P_1\left(\frac{\log x}{\log K}\right)\tilde{f}_1\left(\frac{\log K/(xx_1x_2)}{\log K}\right)
\frac{dx}{x}dx_2\,dx_1+O\left((\log T)^{r^2-1}\right).
\end{split}
\]
By the variable changes $u=1-\frac{\log x_1}{\log K}$, $v=1-\frac{\log x_1x_2}{\log K}$, $w=\frac{\log x}{\log K}$, we get
\be
\label{PreobN43}
\begin{split}
N_{43}=&-\frac{2A_rr^5}{\pi}(\log T)^{r^2}\int_0^1P_1(1-u)\int_0^u(u-v)^{r^2-1}\tilde{f}_1(v)\int_0^v\sin(\pi cw)P_1(w)\tilde{f}_1(v-w)\,dw\,dv\,du\\
&+O\left((\log T)^{r^2-1+\epsilon}\right),
\end{split}
\ee
where the constant in the $O$-term depends on $r$, $\epsilon$ and $\tilde{f}_1$, $P$.

Using $D_i$, $N_i$ given by~\eqref{PreobD1}--\eqref{PreobN43}
we can evaluate
\[
h(c)=c-\frac{N_1+N_2+N_3+N_4}{D_1+D_2+D_3}.
\]
The results of our numerical calculations are summarized in Table~\ref{Preobtab}.

\begin{table}
\begin{center}
\begin{tabular}{|ccc||cc||p{1.2in}|p{1.2in}|c|}\hline
\multicolumn{3}{|c||}{Degrees} & & & \multicolumn{3}{c|}{Polynomials} \\ \cline{1-3}\cline{6-8}
\multicolumn{1}{|c|}{$f_1$} & \multicolumn{1}{c|}{$\tilde{f}_1$} & $P$ & Value of $c$ & Value of $r$ & \multicolumn{1}{c|}{$f_1$} & \multicolumn{1}{c|}{$\tilde{f}_1$} & $P$ \\ \hline\hline
3 & 1 & 2 & $0{.}515398$ & $1{.}18$ & $1{.}95+1{.}47x-1{.}07x^2-0{.}29x^3$ & $-0{.}7-1{.}92x$ & $x^2$ \\ \hline
3 & 1 & 3 & $0{.}515397$ & $1{.}18$ & $1{.}655+1{.}25x-0{.}886x^2-0{.}25x^3$ & $-0{.}57-1{.}6x$ & $x^2+0{.}036x^3$ \\ \hline
6 & 2 & 3 & $0{.}515396$ & $1{.}18$ & $1{.}78+1{.}017x+0{.}2x^2-1{.}56x^3+0{.}45x^4-0{.}06x^5+0{.}05x^6$ & $-0{.}629-0{.}88x-1{.}799x^2$ & $x^2+0{.}083x^3$ \\ \hline
\end{tabular}
\end{center}
\caption{Numerically optimal polynomials in the coefficients $\{a_k\}$, for which $h(c)>1$.}
\label{Preobtab}
\end{table}

\bigskip

\end{document}